\documentclass[10pt]{amsart}
\usepackage{amsmath}
\usepackage{amsxtra}
\usepackage{amscd}
\usepackage{amsthm}
\usepackage{amsfonts}
\usepackage{amssymb}
\usepackage{eucal}

\newtheorem{lem}[subsection]{Lemma}
\newtheorem{prop}[subsection]{Proposition}

\newtheorem{thm}[subsection]{Theorem}

\theoremstyle{definition}

\theoremstyle{remark}

\newcommand{\thmref}[1]{Theorem~\ref{#1}}
\newcommand{\secref}[1]{Sect.~\ref{#1}}
\newcommand{\lemref}[1]{Lemma~\ref{#1}}
\newcommand{\propref}[1]{Proposition~\ref{#1}}

\newcommand{\nc}{\newcommand}
\nc{\renc}{\renewcommand}
\nc{\ssec}{\subsection}
\nc{\sssec}{\subsubsection}
\nc{\on}{\operatorname}

\nc\ol{\overline}
\nc\wt{\widetilde}
\nc\tboxtimes{\wt{\boxtimes}}
\nc{\alp}{\alpha}

\nc{\ZZ}{{\mathbb Z}}
\nc{\NN}{{\mathbb N}}
\nc{\CC}{{\mathbb C}}
\nc{\OO}{{\mathbb O}}
\renc{\SS}{{\mathbb S}}
\nc{\DD}{{\mathbb D}}
\nc{\GG}{{\mathbb G}}

\nc{\Fq}{{\mathbb F}_q}
\nc{\Fqb}{\ol{{\mathbb F}_q}}
\nc{\Ql}{\ol{\mathbb Q}_\ell}
\nc{\id}{\text{id}}
\nc\X{\mathcal X}

\nc{\IC}{\on{IC}}
\nc{\Hom}{\on{Hom}}
\nc{\Lie}{\on{Lie}}
\nc{\Loc}{\on{Loc}}
\nc{\Pic}{\on{Pic}}
\nc{\Bun}{\on{Bun}}
\nc{\ic}{\on{ic}}
\nc{\Aut}{\on{Aut}}
\nc{\rk}{\on{rk}}
\nc{\Sh}{\on{Sh}}
\nc{\Perv}{\on{Perv}}
\nc{\pos}{{\on{pos}}}
\nc{\Conv}{\on{Conv}}
\nc{\Sph}{\on{Sph}}
\nc{\Sym}{\on{Sym}}
\nc{\BunBb}{\overline{\Bun}_B}
\nc{\Buno}{\overset{o}{\Bun}}
\nc{\BunPb}{{\overline{\Bun}_P}}
\nc{\BunBM}{\overline{\Bun}_{B(M)}}
\nc{\BunPbw}{{\widetilde{\Bun}_P}}
\nc{\BunBP}{\widetilde{\Bun}_{B,P}}
\nc{\GUb}{\overline{G/U}}
\nc{\GUPb}{\overline{G/U(P)}}

\nc{\Hhom}{\underline{\on{Hom}}}
\nc\syminfty{\on{Sym}^{\infty}}
\nc\lal{\ol{\lambda}}
\nc\xl{\ol{x}}
\nc\thl{\ol{\theta}}
\nc\nul{\ol{\nu}}
\nc\mul{\ol{\mu}}
\nc\Sum\Sigma
\nc{\oX}{\overset{o}{X}{}}
\nc{\hl}{\overset{\leftarrow}h}
\nc{\hr}{\overset{\rightarrow}h}
\nc{\M}{{\mathcal M}}
\nc{\N}{{\mathcal N}}
\nc{\F}{{\mathcal F}}
\nc{\D}{{\mathcal D}}
\nc{\Q}{{\mathcal Q}}
\nc{\Y}{{\mathcal Y}}
\nc{\G}{{\mathcal G}}
\nc{\E}{{\mathcal E}}
\nc{\CalC}{{\mathcal C}}
\nc\Dh{\widehat{\D}}

\nc{\C}{{\mathcal C}}
\nc{\K}{{\mathcal K}}
\renewcommand{\H}{{\mathcal H}}

\nc{\T}{{\mathcal T}}
\nc{\V}{{\mathcal V}}
\renc{\P}{{\mathcal P}}
\nc{\A}{{\mathcal A}}
\nc{\B}{{\mathcal B}}
\nc{\U}{{\mathcal U}}

\nc{\Gr}{\on{Gr}}

\nc{\frn}{{\check{\mathfrak u}(P)}}
\nc{\p}{\mathfrak p}
\nc{\q}{\mathfrak q}
\nc\f{{\mathfrak f}}

\nc{\qo}{{\mathfrak q}}
\nc{\po}{{\mathfrak p}}
\nc{\s}{{\mathfrak s}}
\nc\w{\text{w}}

\nc\mathi\iota
\nc\Spec{\on{Spec}}
\nc\Mod{\on{Mod}}
\nc{\tw}{\widetilde{\mathfrak t}}
\nc{\pw}{\widetilde{\mathfrak p}}
\nc{\qw}{\widetilde{\mathfrak q}}
\nc{\jw}{\widetilde j}

\nc{\grb}{\overline{\Gr}}
\nc{\I}{\mathcal I}

\nc{\lambdach}{{\check\lambda}}
\nc{\Lambdach}{{\check\Lambda}{}}
\nc{\much}{{\check\mu}}
\nc{\omegach}{{\check\omega}}
\nc{\nuch}{{\check\nu}}
\nc{\etach}{{\check\eta}}
\nc{\alphach}{{\check\alpha}}
\nc{\betach}{{\check\beta}}
\nc{\rhoch}{{\check\rho}}
\nc{\ch}{{\check h}}

\nc{\Hb}{\overline{\H}}

\nc{\BA}{{\mathbb{A}}}
\nc{\BC}{{\mathbb{C}}}
\nc{\BG}{{\mathbb{G}}}
\nc{\BM}{{\mathbb{M}}}
\nc{\BN}{{\mathbb{N}}}
\nc{\BP}{{\mathbb{P}}}
\nc{\BR}{{\mathbb{R}}}
\nc{\BZ}{{\mathbb{Z}}}
\nc{\BS}{{\mathbb{S}}}

\nc{\CA}{{\mathcal{A}}}
\nc{\CB}{{\mathcal{B}}}

\nc{\CE}{{\mathcal{E}}}
\nc{\CF}{{\mathcal{F}}}
\nc{\CG}{{\mathcal{G}}}
\nc{\CL}{{\mathcal{L}}}
\nc{\CM}{{\mathcal{M}}}
\nc{\CN}{{\mathcal{N}}}
\nc{\CK}{{\mathcal{K}}}
\nc{\CO}{{\mathcal{O}}}
\nc{\CP}{{\mathcal{P}}}
\nc{\CQ}{{\mathcal{Q}}}
\nc{\CR}{{\mathcal{R}}}
\nc{\CS}{{\mathcal{S}}}
\nc{\CT}{{\mathcal{T}}}
\nc{\CU}{{\mathcal{U}}}
\nc{\CV}{{\mathcal{V}}}
\nc{\CW}{{\mathcal{W}}}
\nc{\CZ}{{\mathcal{Z}}}
\nc{\CI}{{\mathcal{I}}}

\nc{\cM}{{\check{\mathcal M}}{}}
\nc{\csM}{{\check{\mathcal A}}{}}
\nc{\oM}{{\overset{\circ}{\mathcal M}}{}}
\nc{\obM}{{\overset{\circ}{\mathbf M}}{}}
\nc{\oCA}{{\overset{\circ}{\mathcal A}}{}}
\nc{\obA}{{\overset{\circ}{\mathbf A}}{}}
\nc{\ooM}{{\overset{\circ}{M}}{}}
\nc{\osM}{{\overset{\circ}{\mathsf M}}{}}
\nc{\vM}{{\overset{\bullet}{\mathcal M}}{}}
\nc{\nM}{{\underset{\bullet}{\mathcal M}}{}}
\nc{\oD}{{\overset{\circ}{\mathcal D}}{}}
\nc{\obD}{{\overset{\circ}{\mathbf D}}{}}
\nc{\oA}{{\overset{\circ}{\mathbb A}}{}}
\nc{\op}{{\overset{\bullet}{\mathbf p}}{}}
\nc{\cp}{{\overset{\circ}{\mathbf p}}{}}
\nc{\oU}{{\overset{\bullet}{\mathcal U}}{}}
\nc{\oZ}{{\overset{\circ}{\mathcal Z}}{}}
\nc{\ofZ}{{\overset{\circ}{\mathfrak Z}}{}}
\nc{\oF}{{\overset{\circ}{\fF}}}

\nc{\fa}{{\mathfrak{a}}}
\nc{\fb}{{\mathfrak{b}}}
\nc{\fg}{{\mathfrak{g}}}
\nc{\fgl}{{\mathfrak{gl}}}
\nc{\fh}{{\mathfrak{h}}}
\nc{\fj}{{\mathfrak{j}}}
\nc{\fm}{{\mathfrak{m}}}
\nc{\fn}{{\mathfrak{n}}}
\nc{\fu}{{\mathfrak{u}}}
\nc{\fp}{{\mathfrak{p}}}
\nc{\fr}{{\mathfrak{r}}}
\nc{\fs}{{\mathfrak{s}}}
\nc{\fsl}{{\mathfrak{sl}}}
\nc{\hsl}{{\widehat{\mathfrak{sl}}}}
\nc{\hgl}{{\widehat{\mathfrak{gl}}}}
\nc{\hg}{{\widehat{\mathfrak{g}}}}
\nc{\chg}{{\widehat{\mathfrak{g}}}{}^\vee}
\nc{\hn}{{\widehat{\mathfrak{n}}}}
\nc{\chn}{{\widehat{\mathfrak{n}}}{}^\vee}

\nc{\fA}{{\mathfrak{A}}}
\nc{\fB}{{\mathfrak{B}}}
\nc{\fD}{{\mathfrak{D}}}
\nc{\fE}{{\mathfrak{E}}}
\nc{\fF}{{\mathfrak{F}}}
\nc{\fG}{{\mathfrak{G}}}
\nc{\fK}{{\mathfrak{K}}}
\nc{\fL}{{\mathfrak{L}}}
\nc{\fM}{{\mathfrak{M}}}
\nc{\fN}{{\mathfrak{N}}}
\nc{\fP}{{\mathfrak{P}}}
\nc{\fU}{{\mathfrak{U}}}
\nc{\fV}{{\mathfrak{V}}}
\nc{\fq}{{\mathfrak{q}}}
\nc{\fZ}{{\mathfrak{Z}}}

\nc{\bb}{{\mathbf{b}}}
\nc{\bc}{{\mathbf{c}}}
\nc{\bd}{{\mathbf{d}}}
\nc{\be}{{\mathbf{e}}}
\nc{\bj}{{\mathbf{j}}}
\nc{\bn}{{\mathbf{n}}}
\nc{\bp}{{\mathbf{p}}}
\nc{\bq}{{\mathbf{q}}}
\nc{\bu}{{\mathbf{u}}}
\nc{\bv}{{\mathbf{v}}}
\nc{\bx}{{\mathbf{x}}}
\nc{\bs}{{\mathbf{s}}}
\nc{\by}{{\mathbf{y}}}
\nc{\bw}{{\mathbf{w}}}
\nc{\bA}{{\mathbf{A}}}
\nc{\bK}{{\mathbf{K}}}
\nc{\bB}{{\mathbf{B}}}
\nc{\bC}{{\mathbf{C}}}
\nc{\bG}{{\mathbf{G}}}
\nc{\bD}{{\mathbf{D}}}
\nc{\bH}{{\mathbf{H}}}
\nc{\bM}{{\mathbf{M}}}
\nc{\bN}{{\mathbf{N}}}
\nc{\bV}{{\mathbf{V}}}
\nc{\bW}{{\mathbf{W}}}
\nc{\bX}{{\mathbf{X}}}
\nc{\bZ}{{\mathbf{Z}}}
\nc{\bS}{{\mathbf{S}}}

\nc{\sA}{{\mathsf{A}}}
\nc{\sB}{{\mathsf{B}}}
\nc{\sC}{{\mathsf{C}}}
\nc{\sD}{{\mathsf{D}}}
\nc{\sF}{{\mathsf{F}}}
\nc{\sK}{{\mathsf{K}}}
\nc{\sM}{{\mathsf{M}}}
\nc{\sO}{{\mathsf{O}}}
\nc{\sQ}{{\mathsf{Q}}}
\nc{\sP}{{\mathsf{P}}}
\nc{\sZ}{{\mathsf{Z}}}
\nc{\sfp}{{\mathsf{p}}}
\nc{\sr}{{\mathsf{r}}}
\nc{\sg}{{\mathsf{g}}}
\nc{\sff}{{\mathsf{f}}}
\nc{\sfb}{{\mathsf{b}}}
\nc{\sfc}{{\mathsf{c}}}
\nc{\sd}{{\mathsf{d}}}

\nc{\BK}{{\bar{K}}}

\nc{\tA}{{\widetilde{\mathbf{A}}}}
\nc{\tB}{{\widetilde{\mathcal{B}}}}
\nc{\tg}{{\widetilde{\mathfrak{g}}}}
\nc{\tG}{{\widetilde{G}}}
\nc{\TM}{{\widetilde{\mathbb{M}}}{}}
\nc{\tO}{{\widetilde{\mathsf{O}}}{}}
\nc{\tU}{{\widetilde{\mathfrak{U}}}{}}
\nc{\TZ}{{\tilde{Z}}}
\nc{\tx}{{\tilde{x}}}
\nc{\tbv}{{\tilde{\bv}}}
\nc{\tfP}{{\widetilde{\mathfrak{P}}}{}}
\nc{\tz}{{\tilde{\zeta}}}
\nc{\tmu}{{\tilde{\mu}}}

\nc{\urho}{\underline{\rho}}
\nc{\uB}{\underline{B}}
\nc{\uC}{{\underline{\mathbb{C}}}}
\nc{\ui}{\underline{i}}
\nc{\uj}{\underline{j}}
\nc{\ofP}{{\overline{\mathfrak{P}}}}
\nc{\oB}{{\overline{\mathcal{B}}}}
\nc{\og}{{\overline{\mathfrak{g}}}}
\nc{\oI}{{\overline{I}}}

\nc{\eps}{\varepsilon}
\nc{\hrho}{{\hat{\rho}}}

\nc{\one}{{\mathbf{1}}}
\nc{\two}{{\mathbf{t}}}

\nc{\Rep}{{\mathop{\operatorname{\rm Rep}}}}
\nc{\Tot}{{\mathop{\operatorname{\rm Tot}}}}
\nc{\Ker}{{\mathop{\operatorname{\rm Ker}}}}
\nc{\Hilb}{{\mathop{\operatorname{\rm Hilb}}}}
\nc{\End}{{\mathop{\operatorname{\rm End}}}}
\nc{\Ext}{{\mathop{\operatorname{\rm Ext}}}}
\nc{\CHom}{{\mathop{\operatorname{{\mathcal{H}}\it om}}}}
\nc{\GL}{{\mathop{\operatorname{\rm GL}}}}
\nc{\gr}{{\mathop{\operatorname{\rm gr}}}}
\nc{\Id}{{\mathop{\operatorname{\rm Id}}}}
\nc{\de}{{\mathop{\operatorname{\rm def}}}}
\nc{\length}{{\mathop{\operatorname{\rm length}}}}
\nc{\supp}{{\mathop{\operatorname{\rm supp}}}}

\nc{\Cliff}{{\mathsf{Cliff}}}
\nc{\Fl}{\on{Fl}}
\nc{\Fib}{{\mathsf{Fib}}}
\nc{\Coh}{{\mathsf{Coh}}}
\nc{\FCoh}{{\mathsf{FCoh}}}

\nc{\reg}{{\text{\rm reg}}}

\nc{\cplus}{{\mathbf{C}_+}}
\nc{\cminus}{{\mathbf{C}_-}}
\nc{\cthree}{{\mathbf{C}_*}}
\nc{\Qbar}{{\bar{Q}}}

\nc{\bh}{{\bar{h}}}
\nc{\bOmega}{{\overline{\Omega}}}

\nc{\seq}[1]{\stackrel{#1}{\sim}}

\title{Relation between two geometrically defined bases in representations
of $GL_n$}

\author{A.~Braverman, D.~Gaitsgory and M.~Vybornov}

\dedicatory{To the memory of Iosif Donin}

\email{braval@math.brown.edu, gaitsgde@math.harvard.edu, 
vybornov@math.mit.edu}

\date{October 2004}

\begin{document}

\maketitle

\section{The two bases}

\ssec{Introduction}

Following \cite{Gi} and \cite{BG}, one can realize the irreducible finite-dimensional
representation of $GL_n$, corresponding to a certain Young diagram,
in the top cohomology of (the union of) Springer fibers over a nilpotent 
matrix, whose Jordan decomposition corresponds to this diagram. 
We will review this construction in \secref{Lagrangian} below. In particular, the set
of irreducible components of these Springer fibers provides a basis for this
representation. We call it the Springer basis. \footnote{A clarification is in order:
in \cite{BG} two such geometric realizations were considered. 
The first realization is the one considered \cite{Gi}, and it has been recently 
shown to be connected with the one of Nakajima, cf. \cite{Ma, Sa}. The realization 
considered in the present note is the second one, and it is related to the first
one by the Fourier-Deligne transform on the Lie algebra ${\mathfrak {gl}}_d$.
Therefore, in this paper we do not claim any relation between Nakajima's 
realization of finite-dimensional representations and the one of Mirkovi\'c-Vilonen,
which is discussed below.}

On the other hand, we have the theory of geometric Langlands duality, which realizes
the category of finite-dimensional representations of any reductive group $\check G$
in terms of spherical perverse sheaves on the affine Grassmannian of the Langlands
dual group, $\Gr_G$. In particular, by taking the top (and, in fact, the only non-zero)
cohomology with compact supports
of a given irreducible spherical perverse sheaf $\IC^\lambda$, corresponding to
a dominant coweight $\lambda$, along a semi-infinite orbit $S(\mu)\subset \Gr_G$,
corresponding to a coweight $\mu$, we obtain a vector space, which is canonically 
identified with the weight space $V^\lambda(\mu)$, where $V^\lambda$ is the irreducible
representation of $\check G$ with highest weight $\lambda$. Therefore, the
set of irreducible components of the intersection of $S(\mu)$ with the support of
$\IC^\lambda$, provides a basis for $V^\lambda(\mu)$. We call it the Mirkovi\'c-Vilonen 
(MV for short) basis.

Therefore, it is natural to ask whether for $\check G=GL_n$ (in which is case $G$ is also
$GL_n$), the two bases for $V^\lambda(\mu)$ coincide. The purpose of this note
is to prove this fact. 

\medskip

Let us indicate the strategy of the comparison. First, we interpret each given weight
space $V^\lambda(\mu)$ as a multiplicity space of a given finite-dimensional
representation of $GL_n$ in the tensor product of several other ones. 
This set of multiplicities also acquires a basis via the interpretation of tensor
product in terms of convolution of spherical perverse sheaves on the affine 
Grassmannian. One shows (cf.  \secref{Schur-Weyl conv}) that this basis
tautologically coincides with the Springer basis. 

Thus, we have to relate the two bases, both of which are defined in terms
of the affine Grassmannian. The main idea is to view the two appearences
of the affine Grassmannian separately, and in fact to work on the product of
two copies of $\Gr_{GL_n}$, thought of as $\Gr_{GL_n\times GL_n}$.

Finally, we interpret $GL_n\times GL_n$ as a Levi subgroup of $GL_{2n}$,
and the comparison of bases becomes the corollary of the computation of
the intersection cohomology of the Zastava space, which was carried out in
\cite{BFGM}.

\ssec{Conventions}

We will work with varieties over an algebraically closed field $k$ of characteristic $0$,
but the whole discussion carries over to the characteristic $p$ situation with only
minor modifications. The theory of sheaves that we will consider can be either
holonomic D-modules, or constructible $\ell$-adic sheaves, or, if $k=\BC$, 
constructible sheaves in the analytic topology with characteristic $0$ coefficients.
In terms of notation, we opt for $\ell$-adic sheaves. We do not know what part of the
present discussion carries over to the situation, when coefficients of our sheaves
are torsion.

The methods used in this paper rely substantially on the theory of 
spherical perverse sheaves on the affine Grassmannian. We will not
review this theory here, but rather refer the reader to \cite{MV}. In
addition, we will assume familiarity with the results of the paper \cite{BFGM}.

\ssec{The Spinger basis}   \label{Lagrangian}

Let $V$ be an $n$-dimensional vector space over $\Ql$, and we will 
consider $GL(V)=GL^\vee_n$ as an algebraic group over $\Ql$. We fix
a decomposition of $V$ as a direct sum of $1$-dimensional subspaces
$V=V_1\oplus...\oplus V_n$. This defines a maximal torus 
$\BG_m^{\times n}\simeq T^\vee\subset GL(V)$.

For a non-negative integer $d$ we will be interested in representations
of $GL(V)$ that appear as direct summands in $V^{\otimes d}$.
We will call such representations and their highest weights
$d$-positive. The set of weights $\mu$ of $T^\vee$ that appear in
$V^{\otimes d}$ identifies with the set of $n$-tuples of non-negative
integers $\ol{d}=(d_1,...,d_n)$ such that $\underset{i}\Sigma\, d_i=d$;
it will be denoted by $\Lambda^d$. The subset of dominant coweights
among $\Lambda^d$ will be denoted by $\Lambda^{d,+}$;
explicitly it consists of $n$-tuples $\ol{d}$, satisfying $d_1\geq ...\geq d_n$.
As is well-known, the set $\Lambda^{d,+}$ is 
in bijection with each of the following sets:

\medskip

\noindent (a) Nilpotent conjugacy classes in the Lie algebra
${\mathfrak {gl}}_d$ over $k$ with no more than $n$ Jordan blocks.

\smallskip

\noindent (b) Representations of the symmetric group $\Sigma_d$
over $\Ql$, corresponding to Young diagrams with no more than $n$
rows.

\medskip

For $\lambda\in \Lambda^{d,+}$, we denote by $V^\lambda$ the corresponding
irreducible representation of $GL(V)$, normalized so that its highest weight
line is identified with $\underset{i}\otimes\, V_i^{d_i}$. In particular, for
$\lambda$ of the form $(d,0,...,0)$, the corresponding $V^\lambda$ is the symmetric
power $\Sym^d(V)$. We have 
\begin{equation} \label{SW}
V^\lambda\simeq (\rho^\lambda\otimes V^{\otimes d})^{\Sigma_d},
\end{equation}
where $\rho^\lambda$ is the corresponding irreducible 
representation of $\Sigma_d$.

Let us fix $\lambda\in \Lambda^{d,+}$, and $\mu\in \Lambda^d$,
with $\mu=\ol{d}'=(d'_1,...,d'_n)$. By $V^\lambda(\mu)$ we will denote
the $\mu$-weight space in $V^\lambda$. Let us recall the construction of 
a basis in $V^\lambda(\mu)$, following \cite{BG}.

\medskip

Let $\on{Fl}(\mu)$ denote the partial flag variety, classifying flags
$0=M_0\subset M_1\subset ...\subset M_n=k^d$ in the standard
$d$-dimensional space $k^d$, such that $\dim(M_i/M_{i-1})=d'_i$.
We will denote by $A^\lambda$ an arbitrary element in the nilpotent conjugacy 
class in ${\mathfrak {gl}}_d(k)$ corresponding to $\lambda$.
Let $\on{Fl}(\mu)^{A^\lambda}$ be the Springer fiber over 
$A^\lambda\subset {\mathfrak {gl}}_d$, i.e., the subscheme 
of fixed points of $A^\lambda$ acting naturally on $\on{Fl}(\mu)$.
As is well-known, the dimension of this scheme is 
$d^\lambda=\Sigma\, d_i\cdot i-d$.

According to \cite{BG}, 
$H^{2d^\lambda}(\on{Fl}(\mu)^{A^\lambda},\Ql)\otimes \left(\underset{i}\otimes\, V_i^{d'_i}\right)$
identifies naturally with $V^\lambda(\mu)$. Hence, the set $\bB_{Spr}^\lambda(\mu)$
of irreducible components of  $\on{Fl}(\mu)^{A^\lambda}$ provides a basis for
$V^\lambda(\mu)$, tensored by the inverse of the line $\underset{i}\otimes\, V_i^{d'_i}$. 
(This basis does not depend on the choice of an element
$A^\lambda$ in the corresponding nilpotent orbit, since its centralizer in
$GL_d$ is connected.) 

\ssec{The Springer basis via convolution}   \label{Schur-Weyl conv}

Fix another $n$-dimensional vector space $E$ over $\Ql$, and consider
the corresponding group $GL(E)$. For each $\lambda\in \Lambda^{d,+}$, 
we set $E^\lambda$ to be the irreducible representation of $GL(E)$, defined
as 
$$\Hom_{\Sigma_d}(\rho^\lambda,E^{\otimes d}),$$
where $\rho^\lambda$ is as in \eqref{SW}.

Consider the $\Ql$-vector space
$\on{Sym}^d(V\otimes E)$ as a representation of $GL(V)\times GL(E)$. We have:
\begin{align*}
&\on{Sym}^d(V\otimes E)=\on{Hom}_{\Sigma_d}(\Ql,(V\otimes E)^{\otimes d})\simeq
\on{Hom}_{\Sigma_d\times \Sigma_d}(\on{Ind}^{\Sigma_d\times \Sigma_d}_{\Sigma_d}(\Ql),
(V\otimes E)^{\otimes d})\simeq  \\
&\underset{\lambda\in \Lambda^{d,+}}\oplus\, (V^{\otimes d}\otimes \rho^\lambda)^{\Sigma_d}
\otimes \on{Hom}_{\Sigma_d}(\rho^\lambda,E^{\otimes d})\simeq
\underset{\lambda\in \Lambda^{d,+}}\oplus\, V^\lambda\otimes E^\lambda.
\end{align*}
In particular, we obtain an isomorphism of
$GL(V)$-representations:
\begin{equation} \label{decomp of sym power}
\on{Hom}_{GL(E)}(E^\lambda,\on{Sym}^d(V\otimes E))\simeq V^\lambda.
\end{equation}

\medskip

Fix now a weight $\mu=(d'_1,...,d'_n)$ as above. We have the binomial formula
\begin{equation} \label{binomial}
\on{Hom}_{T^\vee}\left(\underset{i}\otimes\, V_i^{d'_i},
\on{Sym}^d(V\otimes E)\right)\simeq \on{Sym}^{d'_1}(E)\otimes...\otimes 
\on{Sym}^{d'_n}(E).
\end{equation}
Hence, we obtain that
\begin{equation}  \label{V lambda as Hom}
V^\lambda(\mu)\simeq \on{Hom}_{GL(E)}(E^\lambda,
\on{Sym}^{d'_1}(E)\otimes...\otimes \on{Sym}^{d'_n}(E))\otimes \left(\underset{i}\otimes \,
V_i^{d'_i}\right).
\end{equation}

\medskip

Let now $\Gr_E$ be the affine Grassmannian of the group dual to $GL(E)$.
We will think of it as the ind-scheme, classifying lattices $\CM$ in $k^n((t))$. We let
$\Gr_E^-$ denote the "negative" part, i.e., the subscheme, consisting of lattices
that contain the standard lattice $\CM_0=k^n[[t]]$. For a non-negative integer
$d$, we let $\Gr_E^{d,-}$ denote the connected component of $\Gr_E^-$
corresponding to lattices such that $\dim_k(\CM/\CM_0)=d$.

We will denote by $\on{Conv}^k(\Gr_E)$ the $k$-fold convolution diagram
$$\underset{k}{\underbrace{\Gr_E\star ...\star \Gr_E}}\overset{p_k}\longrightarrow \Gr_E,$$
which is in fact isomorphic to the $k$-th direct power of $\Gr_E$.
For a collection of non-negative integers $\ol{d}'=d'_1,...,d'_k$ we will denote by 
$\on{Conv}^{\ol{d'},-}(\Gr_E)$ the subscheme of $\on{Conv}^k(\Gr_E)$, which
classifies $k$-tuples of lattices $(\CM_1,...,\CM_k)$, such that
$\CM_{i-1}\subset \CM_i$ and $\dim_k(\CM_i/\CM_{i-1})=d'_i$.

\medskip

Let $Sph_E$ denote the category of spherical perverse sheaves on
$\Gr_E$. For a dominant weight $\lambda$ of $GL(E)$ we will denote by $\Gr^\lambda_E$ 
(resp., $\ol{\Gr}^\lambda_E$) the 
corresponding Schubert variety (resp., its closure), i.e., the $GL_n[[t]]$-orbit
of the diagonal matrix with entries $(t^{-d_1},...,t^{-d_n})$. It is easy to see that 
$\lambda$ belongs to $\Lambda^{d,+}$ if and only if $\Gr^\lambda_E$ is contained in
$\Gr_E^{d,-}$. Let $\IC^\lambda_E$ denote the corresponding irreducible
object of $Sph_E$. 

For $\lambda$ of the form $(d',0,...,0)$, we have $\ol{\Gr}^\lambda_E=\Gr_E^{d,-}$,
and it is rationally smooth, i.e., the IC sheaf on it is the shifted constant
sheaf, and we will use for it the short-hand notation $\IC^{d'}_E$.

\medskip

Let us recall now the basic fact that the category $Sph_E$ is a tensor category
under convolution, denoted $\F_1,\F_2\mapsto \F_1\star\F_2$,
and as such it is equivalent to the category of finite-dimensional
representations of $GL(E)$, once we identify $E$ with the standard $n$-dimensional
subspace. Under this equivalence, the object $\IC^\lambda_E$ goes over to the irreducible
representation with highest weight $\lambda$, with the trivialized highest weight line.
In particular, for $\lambda\in \Lambda^{d,+}$, the corresponding representation
of $GL(E)$ identifies with $E^\lambda$ above.

Hence, in view of \eqref{V lambda as Hom}, for $\lambda\in \Lambda^{d,+}$,
\begin{equation} \label{V lambda as Hom of sheaves}
V^\lambda(\mu)\simeq \on{Hom}_{Sph_E}(\IC^\lambda_E,
,\IC^{d'_1}_E\star...\star \IC^{d'_n}_E)\otimes \left(\underset{i}\otimes \,
V_i^{d'_i}\right).
\end{equation}

The latter Hom space can be interpreted as follows. Choose a point
$\CM\in  \Gr^\lambda_E$, and consider the scheme 
$$p_n^{-1}(\CM)\cap \on{Conv}^{\ol{d'},-}(\Gr_E)\subset 
\on{Conv}^k(\Gr_E).$$
This scheme is of dimension $d^\lambda$ and its top (=$2d^\lambda$-th)
cohomology identifies canonically with the above Hom.
However, it is easy to see that this scheme identifies also with the 
Springer fiber $\on{Fl}(\mu)^{A^\lambda}$. Indeed, let us denote by 
$\CN^d_n$ the variety of nilpotent $d\times d$-matrices, which have
no more than $n$ Jordan blocks, and let $\widetilde{\CN}^{\ol{d'}}_n\subset
\CN^d_n\times \on{Fl}(\mu)$ be the sub-scheme consisting of pairs
$\left(A\in \CN^d_n,(0=M_0\subset M_1\subset...\subset M_n=k^d)\in \on{Fl}(\mu)\right)$,
such that $A(M_i)\subset M_i$. Then we have a natural smooth and surjective 
map of stacks $\Gr_E^{d,-}\to \CN^d_n/GL_d$, and a Cartesian diagram
$$
\CD
\on{Conv}^{\ol{d'},-}(\Gr_E) @>>> \widetilde{\CN}^{\ol{d'}}_n/GL_d \\
@VVV  @VVV  \\
\Gr_E^{d,-} @>>> \CN^d_n/GL_d.
\endCD
$$

This gives rise to another identification of $V^\lambda(\mu)$ with the cohomology
of the Springer fiber $\on{Fl}(\mu)^{A^\lambda}$. However, by unravelling 
the construction of \cite{BG}, we see that the above identification coincides 
with the one of \secref{Lagrangian}.

\ssec{}

We will now recall the construction of another basis for $V^\lambda(\mu)$
by the so-called MV-cylces. We will consider the affine Grassmannian $\Gr_V$,
corresponding to the group dual to $GL(V)$. We consider the standard
lattice $\CM'_0=k^n[[t]]\subset k^n((t))$, and the positive part $\Gr_V^+$
consists of lattices $\CM'$, such that $\CM'\subset \CM'_0$. Replacing
the subscript $E$ by $V$ we adopt the notations from the previous
subsection, modulo the following difference: for $\lambda\in \Lambda^{d,+}$,
we will denote by $\Gr_V^\lambda$ the $GL_n[[t]]$-orbit of the diagonal
matrix $(t^{d_1},...,t^{d_n})$, and we will denote by $\IC^\lambda_V$
the IC-sheaf on its closure.

Consider now some fixed flag $0=\CM'{}^0_0\subset \CM'{}^1_0\subset...
\subset \CM'{}^n_0=\CM'_0$ with $\CM'{}^0_i/\CM'{}^0_{i-1}$ being a free 
$k[[t]]$-module of rank $1$. For a coweight $\mu=\ol{d}'=(d'_1,...,d'_n)$
we will denote by $S(\mu)\subset \Gr_V$ the corresponding semi-infinite orbit,
i.e., the locally closed subvariety of $\Gr_V$, corresponding to lattices $\CM'$,
such that the induced filtration 
$0=\CM'{}^0\subset \CM'{}^1\subset...
\subset \CM'{}^n=\CM'$, defined by
$\CM'{}^i:=\CM'\cap \CM'{}^i_0$,
is such that $\CM'{}^i/\CM'{}^{i-1}$, viewed as a submodule 
in $\CM'{}^i_0/\CM'{}^{i-1}_0$, equals
$\left(\CM'{}^i_0/\CM'{}^{i-1}_0\right)\cdot t^{d'_i}$.

According to \cite{MV}, Theorem 3.2, 
the intersection $S(\mu)\cap \Gr_V^\lambda$ is of dimension
$$d'{}^\lambda(\mu)=d\cdot (n+1)-\underset{i=1,...,n-1}\Sigma\, (d_i+d'_i)\cdot i.$$
Moreover, by the construction of the fiber functor in Theorem 7.3 of {\it loc.cit.},
$$H_c^{2d'{}^\lambda(\mu)}(S(\mu)\cap \Gr_V^\lambda,\Ql)\simeq
H_c^{d'(\mu)}(S(\mu),\IC_V^\lambda)$$ identifies canonically
with $V^\lambda(\mu)\otimes \left(\underset{i}\otimes \,V_i^{d'_i}\right)^{-1}$,
where $d'(\mu):=d'{}^\mu(\mu)$.

In particular, the set $\bB_{MV}^\lambda(\mu)$ of
irreducible components of $S(\mu)\cap \Gr_V^\lambda$ gives a basis of
$V^\lambda(\mu)\otimes \left(\underset{i}\otimes \,V_i^{d'_i}\right)^{-1}$. The
goal of this note is to prove the following:

\begin{thm}  \label{main}
There exists a natural bijection $\bB_{Spr}^\lambda(\mu)\simeq \bB_{MV}^\lambda(\mu)$,
such that the corresponding basis vectors in
$V^\lambda(\mu)\otimes \left(\underset{i}\otimes \,V_i^{d'_i}\right)^{-1}$ coincide.
\end{thm}

\section{The construction}

\ssec{}

The proof of \thmref{main} will be based on a construction involving
the affine Grassmannian of the group dual to $GL(V)\times GL(E)$.

Consider the product $\Gr^{d,+}_V\times \Gr_E^{d,-}$, and consider
now the scheme, denoted $\CP^d_{loc}$, that classifies the data of triples 
$(\CM',\CM,\alpha)$, where $(\CM',\CM)\in  \Gr^{d,+}_V\times \Gr_E^{d,-}$,
and $\alpha$ is an isomorphism of $k[[t]]$-modules $\CM'_0/\CM'\simeq \CM/\CM_0$.
Let $\pi$ denote the natural projection $\CP^d_{loc}\to \Gr^{d,+}_V\times \Gr_E^{d,-}$,
and $\pi_V$ (resp., $\pi_E$) the further projection onto the $\Gr^{d,+}_V$ (resp., $\Gr^{d,-}_E$) 
factor. 

Note that the map $\pi_V$ is smooth, as $\CP^d_{loc}$ identifies with an open subset
of a vector bundle over $\Gr^{d,+}_V$: indeed, the fiber of $\pi_V$ over a given point
$\CM'\in \Gr^{d,+}_V$ is the set of extensions 
$0\to \CM_0\to \CM\to \CM'_0/\CM'\to 0$, which are torsion-free as $k[[t]]$-modules,
i.e., it embeds into $\on{Ext}^1(\CM'_0/\CM',\CM_0)$.

\medskip

For $\lambda\in \Lambda^{d,+}$ let $\CP^\lambda_{loc}$ denote the
subscheme $\pi_E^{-1}(\Gr_E^\lambda)\subset \CP^d_{loc}$, and for 
a weight $\mu$, let us denote by
$\CP^d_{loc}(\mu)$ the subscheme $\pi_V^{-1}(S(\mu))$. Finally,
let us denote by $\CP^\lambda_{loc}(\mu)$ the intersection
$\CP^\lambda_{loc}\cap \CP^d_{loc}(\mu)$.

\begin{prop}  \label{map factors}
The map $\pi_E$, restricted to $\CP^d_{loc}(\mu)$, factors naturally through
a map $$\pi_E^{\on{Conv}}:\CP^d_{loc}(\mu)\to\on{Conv}^{\ol{d'},-}(\Gr_E),$$
followed by $p_n$, where $\mu=\ol{d}'=(d'_1,...,d'_n)$. 
\end{prop}

\begin{proof}

For $\CM'\in \Gr^{d,+}_V$, let $\CN'$ denote the quotient $\CM'{}^0/\CM'$,
viewed as a $k[[t]]$-module. The flag $0=\CM'{}^0_0\subset \CM'{}^1_0
\subset...\subset \CM'{}^n_0=\CM'_0$ induces on $\CN'$ a flag of $k[[t]]$-modules
$0=\CN'{}^0\subset \CN'{}^1\subset...\subset \CN'{}^n=\CN'$, such that
$\dim_k(\CN'{}^i/\CN'{}^{i-1})=d'_i$. Therefore, for any $(\CM',\CM,\alpha)\in 
\CP^\lambda_{loc}(\mu)$,
we obtain a sequence of lattices
$$\CM_0\subset \CM_1\subset...\subset \CM_n=\CM$$
with $\CM_i/\CM_{i-1}\simeq \CN'{}^i/\CN'{}^{i-1}$, and hence of dimension $d'_i$, i.e.,
we obtain the sought-for point of $\on{Conv}^{\ol{d'},-}(\Gr_E).$

\end{proof}

\begin{lem}  \label{smoothness}
The map  $\pi_E^{\on{Conv}}:\CP^d_{loc}(\mu)\to \on{Conv}^{\ol{d'},-}(\Gr_E)$
constructed above is smooth.  Every fiber  is isomorphic to non-empty
open subset of a vector space of dimension $\underset{i=1,...,n}\Sigma\,
d'_i\cdot (n-i+1)$.
\end{lem}

\begin{proof}

The proof is a corollary of the following observation. Let $\CN'{}^\bullet$ be a flag
of torsion modules over $k[[t]]$, with the $i$-th subquotient of dimension $d'_i$.
Let $\CM'{}^\bullet_0$ be a flag of locally free $k[[t]]$-modules, with the $i$-th
subquotient being locally free of rank $1$. (We can consider the above data
over an arbitrary scheme of paramaters.)

Consider the scheme classifying filtered maps $\CM'{}^\bullet_0\to \CN'{}^\bullet$.
This is a vector scheme of dimension $\underset{i=1,...,n}\Sigma\,
d'_i\cdot (n-i+1)$. The fibers of the map of interest are isomorphic to an 
open subscheme in the above scheme of filtered maps, corresponding to 
the condition that the maps $\CM'{}^i_0/\CM'{}^{i-1}_0\to \CN'{}^i/\CN'{}^{i-1}$
are surjective for all $i=1,...,n$.

\end{proof}

Let us denote by $\bB_{Univ}^\lambda(\mu)$ the set of irreducible components
of the scheme $\CP^\lambda_{loc}(\mu)$. 

\begin{prop} \label{set bijection}
We have canonical bijections
$$\bB_{MV}^\lambda(\mu)\leftarrow 
\bB_{Univ}^\lambda(\mu)\to \bB_{Spr}^\lambda(\mu)$$
\end{prop}

\begin{proof}

The map $\bB_{Univ}^\lambda(\mu)\to \bB_{Spr}^\lambda(\mu)$
comes from the map $\pi_E^{\on{Conv}}$, constructed above.
The fact that it defines a bijection on the corresponding sets 
of irreducible components follows from \lemref{smoothness}.

Let us now construct the map $\bB_{Univ}^\lambda(\mu)\to \bB_{MV}^\lambda(\mu)$.
Note that if $(\CM',\CM,\alpha)\in \CP^d_{loc}$, with $\CM'\in \Gr_V^{\lambda'}$
and $\CM\in \Gr_E^\lambda$, then $\lambda=\lambda'$, because a Schubert cell
contained in $\Gr^{d,+}$ is determined by the isomorphism type of the
quotient $k[[t]]$-module. 

Therefore, the morphism $\pi_V$ restricted to $\CP^\lambda_{loc}$
maps to $\Gr_V^\lambda$ and is a smooth map with non-empty fibers.
Therefore, this map induces a bijection between the set of irreducible 
components of $\CP^\lambda_{loc}(\mu)$ and that of $S(\mu)\cap \Gr_V^\lambda$.

\end{proof}

\ssec{}

The above proposition identifies $\bB_{MV}^\lambda(\mu)$ and $\bB_{Spr}^\lambda(\mu)$
as sets. The rest of the paper is devoted to the identification of the
corresponding basis vectors.

\begin{lem}  \label{dim estimate}
The complex $\pi_!(\Ql)$ is concentrated in the perverse cohomological
degrees $\leq 2\cdot d\cdot n$. 
\end{lem}

The lemma follows from the interpretation of the scheme $\CP^d_{loc}$ as the central
fiber of the Zastava space, explained in the next section, combined with Proposition 5.7
of \cite{BFGM}.

\medskip

Let $Sph_{V\times E}$ denote the category of spherical perverse sheaves
on the affine Grassmannian $\Gr_V\times \Gr_E$. This is a tensor category, equivalent
to the category of representations of the group $GL(V)\times GL(E)$.
The top (=$2\cdot d\cdot n$) cohomology perverse sheaf of the above lemma is clearly spherical;
let us denote it by $\F^d_{loc}$.

For a weight $\mu$ let us denote by $\fq(\mu)$ the projection $S(\mu)\times \Gr_E\to \Gr_E$.
For $\F\in Sph_{V\times E}$ consider its (usual) inverse image to $S(\mu)\times \Gr_E$,
followed by the direct image with compact supports with respect to $\fq(\mu)$,
cohomologically shifted by $[d'(\mu)]$. Slightly abusing the notation, we will
denote the resulting functor $Sph_{V\times E}\to Sph_E$
by $\F\mapsto \fq(\mu)_!(\F)[d'(\mu)]$. (This is indeed a perverse sheaf, i.e., no
lower perverse cohomology appears,  by \cite{MV}, Theorem 3.5.)

In terms of the equivalence $Sph_{V\times E}\simeq \on{Rep}(GL(V)\times GL(E))$,
the above functor corresponds to sending a representation 
to the weight $\mu$ component with respect to the maximal torus
of $GL(V)$, i.e., for $\mu=\ol{d}'=(d'_1,...,d'_n)$ this functor sends a representation $W$
of $GL(V)\times GL(E)$ to the $GL(E)$-representation
$\on{Hom}_{T^\vee}\left(\underset{i}\otimes\, V_i^{d'_i},W\right)$.

\medskip

Using \propref{map factors} we obtain that there exists a canonical isomorphism
\begin{equation} \label{one projection}
\fq(\mu)_!(\F^d_{loc})[d'(\mu)]\simeq \IC^{d'_1}_E\star...\star \IC^{d'_n}_E.
\end{equation}

The following assertion will be proved in the next section:

\begin{thm} \label{identification of sheaf} \hfill

\smallskip

\noindent {\em (a)}
Under the equivalence $\on{Rep}(GL(V)\times GL(E))\simeq Sph_{V\times E}$
the object $\F^d_{loc}$ goes over to $\on{Sym}^d(V\otimes E)$.

\smallskip 

\noindent {\em (b)}
Under this identification, the isomorphism of \eqref{one projection} coincides with
that of \eqref{binomial}.

\end{thm}

\ssec{}

We will now deduce \thmref{main} from \thmref{identification of sheaf}.
By \eqref{decomp of sym power} and
\thmref{identification of sheaf}(a),
\begin{equation} \label{sheaf as direct sum}
\F^d_{loc}\simeq \underset{\lambda\in \Lambda^{d,+}}\oplus\, \IC_V^\lambda\boxtimes \IC_E^\lambda.
\end{equation}

Consider the vector space 
$$\wt{V}{}^\lambda(\mu):=\on{Hom}_{Sph_E}
\left(\IC^\lambda_E,\fq(\mu)_!(\F^d_{loc})[2d'(\mu)])\right).$$
On the one hand, by \eqref{sheaf as direct sum}, we have:
\begin{equation} \label{ident 1}
\wt{V}{}^\lambda(\mu)\simeq H^{d'(\mu)}_c(S(\mu),\IC^\lambda_V)\simeq V^\lambda(\mu)
\otimes \left(\underset{i}\otimes\, V_i^{d'_i}\right)^{-1}.
\end{equation} 
On the other hand, by \eqref{one projection} and \eqref{V lambda as Hom},
\begin{equation} \label{ident 2}
\wt{V}{}^\lambda(\mu)\simeq \on{Hom}_{Sph_E}(\IC^\lambda_E,\IC^{d'_1}_E\star...\star
\IC^{d'_n}_E)\simeq V^\lambda(\mu)\otimes \left(\underset{i}\otimes\, V_i^{d'_i}\right)^{-1}.
\end{equation}
However, by \thmref{identification of sheaf}(b), the resulting two isomorphisms
$\wt{V}{}^\lambda(\mu)\simeq V^\lambda(\mu)\otimes \left(\underset{i}\otimes\, V_i^{d'_i}\right)^{-1}$ coincide.

\medskip

By construction, the set $\bB_{Univ}^\lambda(\mu)$
of irreducible components of $\CP^\lambda_{loc}(\mu)$ defines a basis in 
$\wt{V}{}^\lambda(\mu)$.
We claim that under the identification of \eqref{ident 1}, this basis goes over
to the basis given by $\bB_{MV}^\lambda(\mu)$, and under the identification
of \eqref{ident 2}, $\bB_{Univ}^\lambda(\mu)$ goes over to the basis given by
$\bB_{Spr}^\lambda(\mu)$, where the bijections on the level of underlying sets
are given by \propref{set bijection}. Clearly, this would imply the assertion of the 
\thmref{main}.

\medskip

The first assertion follows readily from the fact that that the map
$$\pi_V:\CP^\lambda_{loc}(\mu)\to \left(S(\mu)\cap \Gr_V^\lambda\right),$$
considered in the proof of \propref{set bijection}, is
a smooth fibration. The second assertion follows from the corresponding
property of the morphism 
$$\pi_E^{\on{Conv}}: \CP^\lambda_{loc}(\mu)\to \left(\on{Conv}^{\ol{d}',-}(\Gr_E)\cap 
p_k^{-1}(\Gr_E^\lambda)\right),$$
given by \lemref{smoothness}.

\section{Proof of \thmref{identification of sheaf} and Zastava spaces}

\ssec{}

Let $G$ be the group $GL_{2n}$, thought of as the dual group of
$GL(V\oplus E)$, let $P$ be the maximal parabolic, whose Levi quotient
$M$ is $GL_n\times GL_n$, thought of as the group dual to $GL(V)\times GL(E)$,
and the unipotent radical $N\simeq \on{Mat}_{n,n}$.
Following \cite{BFGM}, Sect. 2.2,
we will consider the enhanced Zastava space, corresponding
to the pair $(G,P)$. It depends on a paramater $d\in \BZ^+$, and we will denote it 
here by $\ol{\CP}^d_{glob}$. Let us remind the definition in the form adapted to the 
present situation:

Let $X$ be a smooth algebraic curve (not necessarily complete). Let
$\CM_0$ and $\CM'_0$ be the trivial rank $n$-vector bundles on $X$.
The scheme $\ol{\CP}^d_{glob}$ classifies the data of
$$(0\to \CM_0\overset{\imath_0}\to \CA\overset{\jmath_0}\to \CM'_0\to 0,
\CM,\CM',\imath,\jmath),$$
where $0\to \CM_0\to \CA\to \CM'_0\to 0$ is 
a short exact sequence of vector bundles;
$\CM,\CM'$ are rank $n$ vector bundles; $\imath$ and $\jmath$ are maps of
coherent sheaves
$$\imath:\CM'\to \CA \text{ and } \jmath:\CA\to \CM,$$
such that $\jmath\circ \imath=0$, and such that the compositions 
$$\beta:=\jmath\circ \imath_0:\CM_0\to \CM \text{ and }
\beta':=\jmath_0\circ \imath:\CM'\to \CM'_0$$ are both injective. Finally, we require that resulting 
meromorphic isomorphism $\on{det}(\CA)\simeq \on{det}(\CM)\otimes \on{det}(\CM')$
be globally regular.

\medskip

In other words, a point of $\ol{\CP}^d_{glob}$ can be thought of as a diagram
$$
\CD
& &  & & \CM' \\
&  &  & & @V{\imath}VV  \\
0 @>{\imath_0}>> \CM_0 @>>> \CA @>{\jmath_0}>> \CM'_0   @>>>  0 \\
& & & & @V{\jmath}VV \\
& & & & \CM .
\endCD
$$

In particular, we have 
an isomorphism $\CA/(\CM_0\oplus \CM')\simeq \CM'_0/\CM'$ and a map
$\CA/(\CM_0\oplus \CM')\to \CM/\CM_0$. We will denote by $\alpha$ the 
resulting map $\CM'_0/\CM'\to \CM/\CM_0$.

\medskip

Let $\Mod_E^d$ denote the scheme that classifies the data of pairs
$(\CM,\beta)$, where $\CM$ is a rank $n$ vector bundle, and $\beta$ is 
an injective map of coherent sheaves $\CM_0\hookrightarrow \CM$, such that the
quotient torsion sheaf has length $d$. We have a natural map
$\fs_E:\Mod^d_E\to X^{(d)}$.  Similarly, we introduce the scheme $\Mod^d_V$,
which classifies {\it lower} modifications of $\CM'_0$ of length $d$. We obtain
a natural map
$$\ol{\pi}:\ol{\CP}^d_{glob}\to \Mod^d_V\underset{X^{(d)}}\times \Mod_E^d.$$

Note that the data of $\left((\CM,\beta),(\CM',\beta'),\alpha\right)$, such that 
$\fs_E(\CM,\beta)=\fs_V(\CM',\beta')\in X^{(d)}$, is equivalent to
a point of $\ol{\CP}^d_{glob}$. Indeed, one reconstructs $\CA$ as the preimage
under $$\CM'_0\oplus \CM\twoheadrightarrow \CM'_0/\CM'\oplus \CM/\CM_0$$
of the graph of $\CM'_0/\CM'$, embedded by means of $\on{id}\oplus \alpha$.

Let $\CP^d_{glob}$ denote the open subscheme of $\ol{\CP}^d_{glob}$, corresponding
to the condition that the map $\alpha$ is an isomorphism. Note that this is equivalent
to each of the following two conditions: that $\CM'\to \CA$ is a bundle map, or 
that $\CA\to \CM$ is surjective. Therefore, the scheme $\CP^d_{glob}$ classifies the
data of a $G$-bundle on $X$, together with a reduction to $P$ and a reduction to $N$,
which are transversal at the generic point of the curve. We will denote by $\pi$ the
restriction of the morphism $\ol\pi$ to this open subscheme.

\medskip

Let now $x\in X$ be some chosen point, and $t$ a local coordinate. Note that
the preimage of $d\cdot x\in X^{(d)}$ under the map 
$\Mod^d_V\underset{X^{(d)}}\times \Mod_E^d \to X^{(d)}$ identifies
with $\Gr_V^{d,+}\times \Gr_E^{d,-}$; and the preimage of the latter subscheme
under $\CP^d_{glob} \to \Mod^d_V\underset{X^{(d)}}\times \Mod_E^d$
identifies naturally with the scheme $\CP^d_{loc}$, introduced in the previous 
section.

\ssec{}

Let $\IC_{\ol{\CP}^d_{glob}}$ denote the intersection cohomology sheaf
on $\ol{\CP}^d_{glob}$. According to \cite{BFGM}, Cor. 3.8, the open subset  $\CP^d_{glob}$
is smooth, so $\IC_{\CP^d_{glob}}\simeq \Ql[2\cdot d\cdot n]$. Consider the direct image
$$\F^d_{glob}:=\ol{\pi}_!(\IC_{\ol{\CP}^d_{glob}}).$$
The following is the result of application of the machinery of \cite{BFGM}, Theorem 4.5,
Proposition 5.2 and Theorem 5.9 of {\it loc.cit.} to our pair $(G,P)$:

\begin{thm}  \label{Zastava}

\hfill

\smallskip

\noindent{\em (a)}
$\F^d_{glob}$ is isomorphic to the intersection cohomology sheaf
of $\Mod^d_V\underset{X^{(d)}}\times \Mod_E^d$.

\smallskip

\noindent{\em (b)}
The natural map $\pi_!(\Ql[2\cdot d\cdot n])|_{d\cdot x}\to 
\F^d_{glob}|_{d\cdot x}$ induces an isomorphism
$$\F^d_{loc}\simeq \F^d_{glob}|_{d\cdot x}[-d].$$
\end{thm}

\medskip

\thmref{Zastava} readily implies point (a) of \thmref{identification of sheaf}. 
Indeed, let $\overset{\circ}{X}{}^{(d)}$ denote the open subset of $X^{(d)}$, corresponding
to multiplicity-free divisors, and let $\overset{\circ}{X}{}^d$ be its preimage in the
$d$-th Cartesian power of $X$. We have:
$$(\Mod^d_V\underset{X^{(d)}}\times \Mod_E^d)\underset{X^{(d)}}\times
\overset{\circ}{X}{}^d\simeq 
(\Mod^1_V\underset{X}\times \Mod_E^1)^{\times d}\underset{X^d}\times
\overset{\circ}{X}{}^d.$$

For each point $x\in X$, the object of $\on{Rep}(GL(V)\times GL(E))$, corresponding
to 
$$\IC_{\Mod^1_V\underset{X}\times \Mod_E^1}[-1]|_{x}\simeq 
\Ql{}_{\Gr^{1,+}_V\times \Gr^{1,-}_E}[2(n-1)]\in Sph_{V\times E}$$
is $V\otimes E$. Now the assertion of \thmref{identification of sheaf}(a)
follows from Lemma 4.3 of \cite{BFGM}.

\ssec{}

It remains to prove point (b) of \thmref{identification of sheaf}. We fix a flag 
$0=\CM'{}^0_0\subset \CM'{}^1_0\subset...\subset \CM'{}^n_0=\CM'_0$
of sub-bundles in $\CM'_0$, with $\CM'{}^i_0$ being of rank $i$. 
For a weight $\mu=\ol{d}'=(d'_1,...,d'_n)$ let $S_{glob}(\mu)$ be a locally closed
subscheme of $\Mod_V^d$ consisting of pairs 
$(\CM',\beta:\CM'\hookrightarrow \CM'_0)$, such that the induced 
filtration $\CM'{}^i$ on $\CM'$ is such that each subquotient 
$(\CM'{}^i_0/\CM'{}^{i-1}_0)/(\CM'{}^i/\CM'{}^{i-1})$, which is a torsion
sheaf on $X$, is of length $d'_i$. We have a natural map
$\fs(\mu):S_{glob}(\mu)\to X^\mu:=X^{(d'_1)}\times...\times X^{(d'_n)}$. 

By analogy with the local situation, we will denote by
$\F\mapsto \fq(\mu)_!(\F)[d'(\mu)]$ the functor from the derived
category on $\Mod^d_V\underset{X^{(d)}}\times \Mod_E^d$ to that
on $X^\mu\underset{X^{(d)}}\times \Mod_E^d$, given by restriction 
to $S_{glob}(\mu)\underset{X^{(d)}}\times \Mod_E^d$,
followed by the $[d'(\mu)]$-shifted direct image onto 
$X^\mu\underset{X^{(d)}}\times\Mod_E^d$. 

\medskip

Let us also consider the scheme $\on{Conv}^{\ol{d}'}(\Mod_E)$, which classifies the data
of sequences $\CM_0\subset\CM_1\subset...\subset \CM_n=\CM$ of rank $n$ vector
bundles, embedded into one-another as coherent sheaves, such that each
quotient $\CM_i/\CM_{i-1}$ is of length $d'_i$. We will denote by $p_n$
the forgetful map $\on{Conv}^{\ol{d}'}(\Mod_E)\to \Mod_E^d$. By remembering the supports
of the subquotients, we obtain a map $\on{Conv}^{\ol{d}'}(\Mod_E)\to X^\mu$.

Let $\overset{\circ}{X}{}^\mu$ denote the preimage of
this open subset under $X^\mu\to X^{(d)}$. Note that we have natural isomorphisms
\begin{equation} \label{factorization}
\overset{\circ}{X}{}^\mu\underset{X^{(d)}}\times\Mod_E^d\simeq
\underset{i=1,...,n}\Pi\, \overset{\circ}{X}{}^{(d'_i)}\underset{X^{(d'_i)}}\times
\Mod_E^{d'_i}\simeq \overset{\circ}{X}{}^\mu\underset{X^\mu}\times
\on{Conv}^{\ol{d}'}(\Mod_E).
\end{equation}

\medskip

Using \cite{MV}, Theorem 3.6, and \thmref{Zastava}(a) we obtain that 
$\fq(\mu)_!(\F^d_{glob})[d'(\mu)]$ is a perverse sheaf on 
$X^\mu\underset{X^{(d)}}\times\Mod_E^d$, and it equals
the Goresky-MacPherson extension of its restriction to the open subscheme
$\overset{\circ}{X}{}^\mu\underset{X^{(d)}}\times\Mod_E^d$.

The latter restriction identifies, in terms of the isomorphism of \eqref{factorization},
with the constant perverse sheaf on $\overset{\circ}{X}{}^\mu\underset{X^\mu}\times
\on{Conv}^{\ol{d}'}(\Mod_E)$. By the smallness of the map 
$p_n:\on{Conv}^{\ol{d}'}(\Mod_E)\to \Mod_E^d$, we obtain an isomorphism
$$\fq(\mu)_!(\F^d_{glob})[d'(\mu)]\simeq (p_n)_!(\Ql[d\cdot n]).$$ 
By restricting this isomorphism to the 
preimage of $d\cdot x\in X^{(d)}$, we obtain an identification of 
$\fq(\mu)_!(\F^d_{loc})[d'(\mu)]$ with
$\IC^{d'_1}_E\star...\star \IC^{d'_n}_E$, which coincides with that coming from \eqref{binomial},
by the construction of the commutativity constraint on $Sph_E$ via fusion.

\medskip

Let us denote by
$\CP^d_{glob}(\mu)$ the preimage of $S_{glob}(\mu)$ under the natural
projection $\CP^d_{glob}\to \Mod_V^d$.
As in \propref{map factors} and \lemref{smoothness}, the map
$$\pi_E:\CP^d_{glob}(\mu)\to S_{glob}(\mu)\underset{X^{(d)}}\times \Mod_E^d\to
\Mod_E^d$$
factors naturally through a map
$$\pi_E^{\on{Conv}}:\CP^d_{glob}(\mu)\to 
S_{glob}(\mu)\underset{X^\mu}\times \on{Conv}^{\ol{d}'}(\Mod_E)\to \Mod^d_E,$$
followed by $p_n$. 
The top perverse cohomology of $(\pi_E^{\on{Conv}})_!(\Ql)$
is the constant perverse sheaf on $\on{Conv}^{\ol{d}'}(\Mod_E)$. This induces another isomorphism
$$\fq(\mu)_!(\F^d_{glob})[d'(\mu)]\simeq (p_n)_!(\Ql[d\cdot n]).$$
The restriction of this isomorphism 
on the preimage of $d\cdot x\in X^{(d)}$ gives rise to the isomorphism of
\eqref{one projection}.

However, the above two isomorphisms between 
$\fq(\mu)_!(\F^d_{glob})[d'(\mu)]$ and $(p_n)_!(\Ql[d\cdot n])$
coincide, because this is evidently true over the open subscheme
$\overset{\circ}{X}{}^\mu\underset{X^{(d)}}\times\Mod_E^d$. Therefore, their
restrictions to the fibers over $d\cdot x\in X^{(d)}$ coincide too, which is what
we had to show.

\bigskip

\noindent{\bf Acknowledgments.}
We would like to thank Ivan Mirkovi\'c for stimulating discussions.
The work of the first and the second authors is supported by grants
from the NSF. In addition,  the work of the second author is supported 
by a grant from DARPA via NSF, DMS 0105256.

\end{document}